\documentclass[12pt,a4paper]{article}
\usepackage{amsmath,amssymb,amsfonts}
\def\Bbb{\mathbb}

\title{\bf Hecke's inverse cotangent numbers}

\author{Kurt Girstmair}

\date{}

\makeatletter
\let\@@maketitle=\maketitle
\def\maketitle{\def\thispagestyle##1{\relax}\@@maketitle}
\makeatother
\newtheorem{theorem}{Theorem}
\newtheorem{prop}{Proposition}

\newtheorem{corol}{Corollary}

\def\BE{\begin{equation}}
\def\EE{\end{equation}}
\def\BD{\begin{displaymath}}
\def\ED{\end{displaymath}}
\def\BA{\begin{array}}
\def\EA{\end{array}}
\def\BEA{\begin{eqnarray*}}
\def\EEA{\end{eqnarray*}}
\def\BI{\bibitem}

\def\Z{\Bbb Z}
\def\Q{\Bbb Q}

\def\C{\Bbb C}

\def\RR{{\cal R}}
\def\XX{{\cal X}}

\def\phi{\varphi}

\def\MB{\mbox}
\def\LD{\ldots}
\def\OV{\overline}

\def\WH{\widehat}

\def\DIV{\,|\,}

\def\BQ{``}
\def\EQ{'' }
\def\EQP{''}

\def\MN{\medskip\noindent}
\def\STOP{\hfill$\Box$}

\def\DR{d^{(r)}}
\def\ctr{{\rm{ct}}^{(r)}}
\def\ctrh{\WH{\rm{ct}}^{(r)}}
\def\CTR{{\rm{Ct}}^{(r)}}
\def\CTRH{\WH{\rm{Ct}}^{(r)}}

\def\BR{B^{(r)}}
\def\BRT{\widetilde B^{(r)}}

\def\NT{\widetilde n}
\def\BRTH{\widehat B^{(r)}}

\def\sh{\widehat s}
\def\ch{\widehat c}

\def\sa{s^{*}}
\def\ca{c^{*}}
\def\ua{\underline a}

\def\XXR{\XX^{(r)}}

\def\BRCHIF{B^{(r)}_{\chi_f}}
\def\BRCHIFB{B^{(r)}_{\OV{\chi}_f}}
\def\B3{B^{(3)}}
\def\ct3h{\WH{\rm{ct}}^{(3)}}

\newcommand{\btop}[2]{\genfrac{}{}{0pt}{1}{#1}{#2}}

\begin{document}
\maketitle

\begin{abstract}

\noindent
In connection with Eisenstein series for the principal congruence subgroup $\Gamma(n)$, Hecke introduced certain numbers, of which he said that they are
rational and cumbersome to calculate. We show, however, that these numbers are (essentially) generators of the $n$th cyclotomic field or its maximal real
subfield. They arise from the well investigated {\em cotangent numbers} by matrix inversion, which is why we call them {\em inverse cotangent numbers}.
We describe them as linear combinations of roots of unity with rational coefficients in a fairly closed form, provided that $n$ is square-free.
We also exhibit a formula for these numbers in terms of generalized Bernoulli numbers and Gauss sums.
\end{abstract}

%%%%%%%%%%%%%%%%%%%%%%%%%%%%%%%%%%%%%%%%%%%%%
\section*{1. Introduction}
%%%%%%%%%%%%%%%%%%%%%%%%%%%%%%%%%%%%%%%%%%%%%

Let $n\ge 3, r\ge 1$ be positive integers. Let $j\in \Z$ be prime to $n$.
In his paper \cite{He} Hecke introduced the numbers
\BE
\label{1.2}
  \DR_j=\sum_{\btop{m\in\Z}{m\equiv j\,{\rm mod }\,n}}\frac{\mu(|m|)}{m^r},
\EE
where $\mu$ denotes the M\"obius function.
In the case $r=1$ the summation of the series has to be understood in the sense of
\BE
\label{1.3}
  \lim_{k\to \infty}\sum_{|m|\le k}.
\EE

Hecke said that the numbers $\DR_j$ turn out to be rational (\BQ ergeben sich als rationale Zahlen\EQP).
However, he doubtlessly knew that a certain power of $\pi$ is involved in $\DR_j$, which, by the analogy of $1/\zeta(r)$, can only be $\pi^{-r}$. Indeed, put
\BE
\label{1.4}
  \ctrh_j=\frac{-i^r\pi^r}{(r-1)!\,n^r}\,\DR_j.
\EE
Then this number generates the $n$th cycotomic field (over $\Q$) if $r$ is odd, and its maximal real subfield if $r$ is even (see Proposition \ref{p1.1}).
In particular, $\ctrh$ is rational only in few cases like $n\in \{3,4,6\}$, $r$ even.
The notation $\widehat{\rm{ct}}$ will turn out to be reasonable.

Hecke further said that the computation of
$\DR_j$ is fairly laborious if $n$ is not a prime, and will not be discussed in his paper (\BQ Ihre Berechnung ist, wenn nicht gerade $N$ eine Primzahl ist,
ziemlich umst\"andlich, und soll hier nicht weiter er\"ortert werden\EQ (Hecke's $N$ corresponds to our $n$)).
We will see, however, that even in the case of a small prime number $n$ and a small parameter $r$ (like $n=11$, $r=4$) the result is fairly complex (see Section 4).

Let us briefly sketch Hecke's context. He studies the vector space of Eisenstein series of weight $r$, $r\in\Z$, $r\ge 1$, for the principal congruence subgroup
$\Gamma (n)=\{A\in \MB{SL}(2,\Z); A\equiv I \MB{ mod } n\}$, where $I$ denotes the unit matrix. Since the cases $r=1,2$ are more complicated, we assume $r\ge 3$ for our sketch.
Hecke's primitive Eisenstein series $G(z,r,\ua)$,
$\ua=(a_1,a_2)$ a pair of coprime integers, have an easily computable Fourier expansion, but they do not form a $\C$-basis of the vector space in question. Therefore, he defines
\BD
   G^*(z,r,\ua)=\sum_{\btop{1\le j\le n/2}{(j,n)=1}}\DR_{j^*}G(z,r,j\ua),
\ED
where $j^*$ is an inverse of $j$ mod $n$,
and obtains such a basis, provided that the pairs $\ua$ are suitably chosen with respect to cusps.

Since the numbers $\ctrh_j$ are in the $n$th-cyclotomic field $\Q(\zeta_n)$, $\zeta_n=e^{2\pi i/n}$, a computation of these number means a representation as a rational linear combination of powers of $\zeta_n$.
We obtain this, in a fairly explicit form, only if $\zeta_n$ generates a {\em normal basis} of $\Q(\zeta_n)$. This, however, is true if, an only if, $n$ is {\em square-free} (see \cite{Leo}).
\BQ Fairly explicit\EQ means up to the inversion of an easily computable rational matrix involving values of Bernoulli polynomials (see Theorem \ref{t1} and the examples of Section 4).

The numbers $\ctrh_j$ are closely connected with the {\em cotangent numbers} $\ctr_j$,
defined by
\BE
\label{1.6}
 \ctr_j=i^r\cot^{(r-1)}(\pi j/n),\enspace (j,n)=1,
\EE
where $\cot^{(r-1)}$ stands for the $(r-1)$th derivative of the cotangent function.
So we have $\cot^{(0)}=\cot$, $\cot^{(1)}=-\cot^2-1$, and so on.
Since the numbers $\ctrh_j$ arise from the numbers $\ctr_j$ by matrix inversion,
they will be called {\em inverse cotangent numbers} (see (\ref{2.8})).

Cotangent numbers have been studied in the case $r=1$ in connection with
formulas for the {\em relative class number} of $\Q(\zeta_n)$, see, for instance, \cite[formula (3.1)]{CaOl}, \cite[formulas (3), (4)]{Gi3}, \cite{En}, \cite{Le}, \cite{IsKa}, for $r\ge 2$ also
\cite{Gi2} and \cite{Gi3}. See also the final remark of this paper.

As a consequence of our study, we exhibit the values of series like
\BD
 \sum_{\btop{m>0,}{(m,n)=1}}\mu(m)\frac{\sin(2\pi m^*/n)}{m^r}, \sum_{\btop{m>0,}{(m,n)=1}}\mu(m)\frac{\cos(2\pi m^*/n)}{m^r},
\ED
where $n$ is square-free, $r$ odd in the first, $r$ even in the second case ($m^*$ has been defined above). These series have rational values up to the factor $\pi^{-r}$, see Proposition \ref{p4}.

Finally, we determine the eigenvalues of the matrices used for the definition of the numbers $\ctrh_j$, see Section 6. This gives another formula for the numbers $\ctrh_j$ in terms of generalized Bernoulli numbers
and Gauss sums, see Theorem \ref{t2}.

%%%%%%%%%%%%%%%%%%%%%%%%%%%%%%%%%%%%%%%%%%%%%
\section*{2. Inverse cotangent numbers}
%%%%%%%%%%%%%%%%%%%%%%%%%%%%%%%%%%%%%%%%%%%%%

Let $n, r$ be as above.
The partial fraction decomposition of the cotangent function yields, for $j\in\Z$, $(j,n)=1$,
\BD
%\label{2.2}
 \cot^{(r-1)}(\pi j/n)=(-1)^{r-1}(r-1)!\left(\frac{n}{\pi}\right)^r\sum_{\btop{m\in\Z}{m\equiv j\,{\rm mod }\,n}}\frac{1}{m^r}.
\ED
In the case $r=1$ the series has to be understood in the sense of (\ref{1.3}).

\begin{prop} % Proposition 1 %%%%%%%%%%%%%%%%%%%%%%%%%%%%
\label{p1}

Let $n$ and $r$ be as above.
Let $\RR=\{j; 1\le j\le n/2, (j,n)=1\}.$ For a number $j\in\Z$, $(j,n)=1$, let $j^*$ denote the inverse of $j$ mod $n$, i.e.,
$jj^*\equiv 1\mod n$, $1\le j^*\le n$.
If $\ctr_j$ and $\ctrh_j$ are defined by \rm{(\ref{1.2}), (\ref{1.4})}, and \rm{(\ref{1.6})}, then
\BE
\label{2.4}
  \sum_{l\in\RR}\ctr_{sl^*}\,\ctrh_{lt^*}=\delta_{s,t},
\EE
for all $s,t\in\RR$. Here $\delta$ is the Kronecker delta.

\end{prop} %%%%%%%%%%%%%%%%%%%%%%%%%%%%%%%%%%%%%%%%%%%%%

\MN
{\em Proof.} The case $r=1$ of conditional convergence is somewhat subtle. It was settled in \cite{Gi1}. In the case $r\ge 2$ the series for $\ctr_j$, $\ctrh_k$ are absolutely
convergent and can, thus, be multiplied without problems. We obtain that the sum in question equals
\BD
  \sum_{\btop{m\in\Z}{m\equiv st^*\,{\rm mod }\,n}}\frac{1}{m^r}\sum_{l\in\RR}\sum_{\btop{j\in\Z,j\DIV m}{j\equiv lt^*\,{\rm mod }\,n}}\mu(|j|).
\ED
For a fixed number $m$, $m\equiv st^*\mod n$, the inner double sum can be written
\BD
  \sum_{l\in\RR}\sum_{\btop{j\DIV m,j>0}{j\equiv lt^*\,\rm{mod}\,n}}\mu(j)+\sum_{l\in\RR}\sum_{\btop{j\DIV m,j>0}{j\equiv -lt^*\,\rm{mod}\,n}}\mu(j).
\ED
This is the same as
\BD
  \sum_{j\DIV m,j>0}\mu(j)=\left\{
                                   \begin{array}{ll}
                                     1, & \hbox{ if $m=\pm 1$;} \\
                                     0, & \hbox{otherwise.}
                                   \end{array}
                                 \right.
\ED
However, the case $m=-1$ would imply $s\equiv -t$ mod $n$, which is impossible for $s,t\in\RR$. Therefore, only $m=1$ and $s\equiv t$ mod $n$ contribute to the whole.
\STOP

\MN
In view of (\ref{2.4}), we introduce the matrices
\BD
\label{2.6}
 \CTR=(\ctr_{jk^*})_{j,k\in\RR},\enspace \CTRH=(\ctrh_{jk^*})_{j,k\in\RR},
\ED
the rows and columns of these matrices having the natural order of $\RR$.
Then Proposition \ref{p1} can be phrased as
\BE
\label{2.8}
  \CTRH=(\CTR)^{-1}.
\EE
This justifies  the name {\em inverse cotangent numbers} given to the numbers $\ctrh_j$ in Section 1.

We further note that the numbers $\ctr_j$, and, by (\ref{2.8}), also the numbers $\ctrh_j$, $(j,n)=1$, lie in $\Q(\zeta_n)$.
Indeed, $\MB{ct}^{(1)}_j=i\cot(\pi j/n)=(1+\zeta_n^j)/(1-\zeta_n^j)$, and $\ctr_j$ is a polynomial in $i\cot(\pi j/n)$ with rational coefficients.

We also need the action of the Galois group of $\Q(\zeta_n)$ over $\Q$ on the numbers $\ctr_j$, $(j,n)=1$.
Let $\sigma_k$, $(k,n)=1$, be the Galois automorphism of $\Q(\zeta_n)$ defined by $\zeta_n\mapsto\zeta_n^k$.
Then
\BE
  \label{2.9}
  \sigma_k(\ctr_j)=\ctr_{kj},\enspace (j,n)=1,
\EE
The action of $\sigma_k$ on the numbers $\ctrh_j$, $j\in\RR$, can be found by (\ref{2.8}), namely, by the identity
\BD
\sigma_k(\CTRH)\cdot \sigma_k(\CTR)=\sigma_k(\CTRH)\cdot(\ctr_{jkl^*})_{j,l\in\RR}=I,
\ED
where $I=(\delta_{j,l})_{j,l\in \RR}$ is the unit matrix. We obtain
\BE
\label{2.9.1}
 \sigma_{k}(\ctrh_j)=\ctrh_{jk^*}.
\EE
The numbers $\ctr_j$, $j\in\RR$, are $\Q$-linearly independent. Indeed, a nontrivial relation
\BD
  \sum_{j\in\RR}b_j\ctr_j=0, \enspace b_j\in\Q,
\ED
would, by Galois action, entail the relations
\BD
  \sum_{j\in\RR}b_j\ctr_{jk^*}=0, \enspace k\in\RR.
\ED
Accordingly, the rows of the matrix $\CTR$ would be linearly dependent (over $\C$), which is impossible,
since the matrix $\CTR$ is invertible. In the same way, the numbers $\ctrh_j$, $j\in \RR$, are $\Q$-linearly independent.

Further, we observe that $\sigma_{-1}\ctr_j$=$(-1)^r\ctr_j$, $j\in\RR$.  In view of (\ref{2.8}), this implies
$\sigma_{-1}(\ctrh_j)=(-1)^r\ctrh_j$. Hence the numbers $\ctrh$ are real if $r$ is even, and purely imaginary if $r$ is odd.

\begin{prop} % Proposition 1.1 %%%%%%%%%%%%%%%%%%%%%%%%%%%%
\label{p1.1}

We adopt the above notation.
Then for each $j$ with $(j,n)=1$,
\BD
   \Q(\ctrh_j)=\left\{
                             \begin{array}{ll}
                               \Q(\zeta_n), & \hbox{if $r$ is odd;} \\
                               \Q(\zeta_n+\zeta_n^{-1}), & \hbox{if $r$ is even.}
                             \end{array}
                           \right.
\ED
\end{prop} %%%%%%%%%%%%%%%%%%%%%%%%%%%%%%%%%%%%%%%%%%%%%

\MN
{\em Proof.}
We know that the numbers $\ctrh_j$, $j\in\RR$, span a $\Q$-vector space of dimension $|\RR|=\phi(n)/2$.
If $r$ is odd, this vector space is not a subfield of $\Q(\zeta_n)$ because all of its numbers are purely imaginary. Since the field
$\Q(\ctrh_j)$ is normal over $\Q$, it contains this subspace, and, accordingly, has a $\Q$-dimension $>\phi(n)/2$.
Therefore, this field coincides with $\Q(\zeta_n)$. If $r$ is even, this subspace lies in the maximal real subfield $\Q(\zeta_n+\zeta_n^{-1})$ and has
the same $\Q$-dimension.
\STOP

Our plan for the next section is as follows. We assume that $n$ is square-free and express the entries of $\CTR$ as rational linear combinations of the basis vectors $\zeta_n^k$, $1\le k\le n$, $(k,n)=1$,
of the field $\Q(\zeta_n)$. This will give the respective representation of the entries of $\CTRH$.

%%%%%%%%%%%%%%%%%%%%%%%%%%%%%%%%%%%%%%%%%%%%%
\section*{3. Generalized Bernoulli matrices}
%%%%%%%%%%%%%%%%%%%%%%%%%%%%%%%%%%%%%%%%%%%%%

We need the following result, which holds for all $n\ge 3$ and seems to be not widely known (see \cite[p. 381]{Gi2} and the reference given there).
\BE
\label{3.2}
 \ctr_1=\frac{2^rn^{r-1}}{r}\sum_{k=1}^n\BR(k/n)\zeta_n^{-k}-\delta_{r,1},
\EE
where $\BR$ denotes the $r$th {\em Bernoulli polynomial}. Here the cotangent number $\ctr_j$ is expressed as a rational linear combination of powers of $\zeta_n$, which are, in general, not $\Q$-linearly independent.

\begin{prop} % Proposition 2 %%%%%%%%%%%%%%%%%%%%%%%%%%%
\label{p2}
If $n\ge 3$ is square-free, the number $\ctr_1$ can be written
\BE
\label{3.4}
  \ctr_1=\sum_{\btop{1\le j\le n}{(j,n)=1}}\BRT_j\zeta_n^j,
\EE
where $\BRT_j$ is given by
\BD
%\label{3.6}
 \BRT_j=\frac{(-1)^r2^rn^{r-1}}{r}\sum_{d\DIV n}\mu(d)\sum_{\btop{1\le k\le \NT}{(k,n)=d,\, k\equiv j\,\rm{mod}\, n/d}}\BR(k/n).
\ED
Here
\BD
  \NT=\left\{
        \begin{array}{ll}
          n-1, & \hbox{if $r$ is odd;} \\
          n, & \hbox{if $r$ is even.}
        \end{array}
      \right.
\ED
\end{prop} %%%%%%%%%%%%%%%%%%%%%%%%%%%%%%%%%%%%%%%%%%%%%

\MN
{\em Proof.}
Using the identity $\BR(1-x)=(-1)^r\BR(x)$, we may write (\ref{3.2}) as
\BD
%\label{3.8}
 \ctr_1=\frac{(-1)^r2^rn^{r-1}}{r}\sum_{k=0}^{n-1}\BR(k/n)\zeta_n^{k}-\delta_{r,1}.
\ED
If $r=1$, $(-1)^r2^rn^{r-1}\BR(0)/r-\delta_{r,1}=0$. If $r$ is odd, $r\ge 3$, $\BR(0)=0$.
Hence we have to deal with the sum
\BD
  \sum_{k=1}^{\NT}\BR(k/n)\zeta_n^{k}
\ED
(observe $\BR(0)=\BR(1)$ if $r$ is even).
This sum can be written
\BE
\label{3.10}
\sum_{d\DIV n}\sum_{\btop{1\le k\le \NT}{(k,n)=d}}\BR(k/n)\zeta_{n/d}^{k/d}.
\EE
Since $n$ is square-free, $\zeta_{n/d}^{k/d}$ is - up to the factor $\mu(d)$ -- the conjugate of a certain trace of $\zeta_n$ (see equation (34) in \cite{Gi3}). Therefore,
\BD
\label{3.12}
  \zeta_{n/d}^{k/d}=\mu(d)\sum_{\btop{1\le j\le n}{(j,n)=1, j\equiv k\,\rm{mod}\,n/d}}\zeta_n^j.
\ED
Accordingly, the sum of (\ref{3.10}) equals
\BD
\label{3.14}
 \sum_{\btop{1\le j\le n}{(j,n)=1}}\zeta^j \sum_{d\DIV n}\mu(d)\sum_{\btop{1\le k\le \NT}{(k,n)=d,  k\equiv j\,\rm{mod}\, n/d}}\BR(k/n),
\ED
which gives the desired result.
\STOP

We define, for $j\in\Z$, $(j,n)=1$,
\BD
\label{3.16}
  s_j=\zeta_n^j-\zeta_n^{-j}=2i\sin(2\pi j/n) \MB{ and } c_j=\zeta_n^j+\zeta_n^{-j}=2\cos(2\pi j/n).
\ED
Recall the definition of $\RR$ from the foregoing section.
Due to the identity $\BR(1-x)=(-1)^r\BR(x)$, we have $\BRT_{n-j}=(-1)^r\BRT_j$.
Hence (\ref{3.4}) yields the following corollary.

\begin{corol} % Corol 1 %%%%%%%%%%%%%%%%%%%%%%%%%%%
\label{c1}
If $n\ge 3$ is square-free, the number $\ctr_1$ can be written
\BE
\label{3.18}
  \ctr_1=\sum_{j\in\RR}\BRT_js_j, \MB{ if } r \MB{ is odd},
\EE
and
\BE
\label{3.20}
  \ctr_1=\sum_{j\in\RR}\BRT_jc_j, \MB{ if } r \MB{ is even}.
\EE

\end{corol} %%%%%%%%%%%%%%%%%%%%%%%%%%%%%%%%%%%%%%%%%%%%%

\MN
Next we introduce certain additional matrices.
To this end we observe
\BD
 \sigma_k(s_j)=s_{kj}, \sigma_k(c_j)=c_{kj},\enspace j\in \Z, (j,n)=1,
\ED
where $\sigma_k$ is the Galois automorphism of the foregoing section.
In view of (\ref{2.9}), (\ref{3.18}), and (\ref{3.20}), we obtain
\BD
  \ctr_k=\left\{
           \begin{array}{ll}
            \sum_{j\in\RR}\BRT_js_{kj} , & \hbox{if $r$ is odd;} \\
             \sum_{j\in\RR}\BRT_jc_{kj}, & \hbox{if $r$ is even.}
           \end{array}
         \right.
\ED
This can be written
\BE
\label{3.23}
  \ctr_k=\sum_{j\in\RR}\BRT_{jk^*}s_{j},\enspace \ctr_k=\sum_{j\in\RR}\BRT_{jk^*}c_{j},
\EE
in the respective cases. Here, however, we have to
make the following {\em convention}: The subscript $jk^*$ means the representative of $jk^*$ in $\{ 1,\LD,n\}$, i.e.,
the number $l$ in this range, defined by $l\equiv jk^*$ mod $n$.

Now we define the matrices
\BD
  \BRT=(\BRT_{jk^*})_{j,k\in\RR},\enspace S=(s_{jk^*})_{j,k\in\RR}, \MB{ and } C=(c_{jk^*})_{j,k\in\RR}.
\ED
We call $\BRT$ the {\em generalized Bernoulli matrix}.
The identities of (\ref{3.23}) yield the following proposition.

\begin{prop} % Proposition 3 %%%%%%%%%%%%%%%%%%%%%%%%%%%
\label{p3}

If $n\ge 3$ is square-free, we have
\BD
%\label{3.24}
  \CTR=\left\{
         \begin{array}{ll}
           S\cdot(\BRT)^T, & \hbox{if $r$ is odd;} \\
           C\cdot(\BRT)^T, & \hbox{if $r$ is even,}
         \end{array}
       \right.,
\ED
where $(\enspace)^T$ denotes the transpose of the respective matrix.

\end{prop} %%%%%%%%%%%%%%%%%%%%%%%%%%%%%%%%%%%%%%%%%%%%%

\MN
Since $\CTRH$ is the inverse of the matrix $\CTR$, we immediately obtain, if $n$ is square-free,
\BE
\label{3.26}
  \CTRH=\left\{
         \begin{array}{ll}
           (\BRTH)^T\cdot S^{-1}, & \hbox{if $r$ is odd;} \\
           (\BRTH)^T\cdot C^{-1}, & \hbox{if $r$ is even,}
         \end{array}
       \right.,
\EE
where $\BRTH=(\BRT)^{-1}$. We write $\BRTH_{j,k}$, $j,k\in\RR$, for the entries of the matrix $\BRTH$.

The inverse matrices of $S$ and $C$ have been given, in the square-free case, in \cite{Gi4}. For the convenience of the reader we repeat their definition here.
For $k\in\Z$, $(k, n)=1$, put
\BD
  \lambda(k)=|\{ q; q\ge 3, q\DIV n, k\equiv 1 \MB{ mod } q \}|.
\ED
Moreover, define, for $k\in\Z$, $(k, n)=1$,
\BD
 \sh_k=\frac{-1}{n}\sum_{l\in\RR}(\lambda(kl)-\lambda(-kl))s_l
\ED
and
\BD
 \ch_k=\frac{1}{n}\sum_{l\in\RR}(\lambda(kl)+\lambda(-kl)+\rho_n)c_l,
\ED
where
\BD
 \rho_n=\left\{
          \begin{array}{ll}
            2, & \hbox{if $n$ is odd;} \\
            4, & \hbox{if $n$ is even.}
          \end{array}
        \right.
\ED
Then
\BD
%\label{3.28}
 S^{-1}=(\sh_{jk^*})_{j,k\in\RR},\enspace C^{-1}=(\ch_{jk^*})_{j,k\in\RR}.
\ED
Now we are in a position to express the numbers $\ctrh_j$, $j\in \RR$, in an explicit way, up to the inversion of the matrx $\BRT$.

\begin{theorem} % Theorem 1 %%%%%%%%%%%%%%%%%%%%%%%%%%%
\label{t1}

If $n\ge 3$ is square-free, we have
\BD
 \ctrh_1=\left\{
           \begin{array}{ll}
             \sum_{j\in\RR}\BRTH_{j,1}\sh_j, & \hbox{ if $r$ is odd;} \\
             \sum_{j\in\RR}\BRTH_{j,1}\ch_j, & \hbox{ if $r$ is even.}
           \end{array}
         \right.
\ED
\end{theorem}

\MN
Theorem \ref{t1} allows writing $\ctrh_1$ in the form
\BD
  \ctrh_1=\sum_{j\in\RR}b_js_j,\: \MB{ and }\:  \ctrh_1=\sum_{j\in\RR}b_jc_j,\enspace b_j\in\Q,
\ED
if $r$ is odd  and $r$ is even, respectively. If we apply $\sigma_{k^*}$, $(k,n)=1$, we obtain
\BD
  \ctrh_k=\sum_{j\in\RR}b_js_{jk^*},\enspace  \ctrh_k=\sum_{j\in\RR}b_jc_{jk^*}, \MB{ respectively},
\ED
(recall (\ref{2.9.1})).

In other words, the number $\ctrh_k$ has, as a linear combination of $s_j$ and $c_j$, $j\in\RR$,
(up to sign changes, if $r$ is odd) the same coefficients as $\ctrh_1$, but in a permuted order.
Hence it suffices to render $\ctrh_1$ in this form, as in Theorem \ref{t1} and in the examples of the next section.

\MN
{\em Remark.} Let $n=p\ge 3$ be a prime. Then the entry $\BRT_{jk^*}$ of the matrix $\BRT$ simplifies to
\BD
  (-1)^r\frac{2^rp^{r-1}}{r}\BR(jk^*/p)\: \MB{ and }\: (-1)^r\frac{2^rp^{r-1}}{r}(\BR(jk^*/p)+\BR(1))
\ED
if $r$ is odd and even, respectively (recall the convention for $jk^*$ made in connection with (\ref{3.23})). Therefore, the matrix $\BRT$ is
basically the matrix of a generalized {\em Maillet determinant} (see, for instance, \cite{CaOl} and \cite{KaKu}).

On the other hand, the numbers $\sh_j$ and $\ch_j$ simplify to
\BE
\label{3.30}
   \sh_j=-s_{j^*}/p \:\MB{ and }\: \ch_j=(c_{j^*}-2)/p,\: (j,p)=1.
\EE
In view of Theorem \ref{t1}, this means that $\ctrh_1$ can be found on combining certain entries of the inverse of a simple Maillet-type matrix
with the numbers of (\ref{3.30}).

%%%%%%%%%%%%%%%%%%%%%%%%%%%%%%%%%%%%%%%%%%%%%
\section*{4. Some examples}
%%%%%%%%%%%%%%%%%%%%%%%%%%%%%%%%%%%%%%%%%%%%%

In this section we give a list of the numbers $\ctrh_1$ for $n$ square-free, $11\le n\le 15$, and $1\le r\le 4$ (up to one exception for reasons of magnitude).
It is more convenient to work with $\sa_j=i\sin(\pi j/n)$ and $\ca_j=\cos(\pi j/n)$ instead of $s_j$ and $c_j$. The cases $n=11, 13$, $r=4$, make it hard to believe that
their results could be found quickly in Hecke's time.

$n=11$: \rule{0mm}{5mm}

$r=1$:  \rule{0mm}{5mm} $- \sa_1/11 -\sa_2/11  - \sa_4/11- \sa_5/11 $

$r=2$: \rule{0mm}{5mm} $ 13\ca_1/275+9\ca_2/550 + 21\ca_3/550 - 3\ca_4/275 + \ca_5/50$

$r=3$: \rule{0mm}{5mm} $79\sa_1/42933+ 409\sa_2/171732 + 41\sa_3/15612 +  161\sa_4/42933 + 769\sa_5/171732$

$r=4$: \rule{0mm}{5mm} $ - 431305\ca_1/568881181 + 27303\ca_2/4551049448 - 2923827\ca_3/4551049448 $ \\ \indent \rule{0mm}{5mm} $+ 256209\ca_4/1137762362 - 1899515\ca_5/4551049448$

$n=13$:\rule{0mm}{5mm}

$r=1$: \rule{0mm}{5mm} $-\sa_2/13 - \sa_3/13  - \sa_5/13- \sa_6/13 $

$r=2$: \rule{0mm}{5mm} $ 103\ca_1/3458 + 37\ca_2/3458  + 75\ca_3/3458- 6\ca_4/1729 + 34\ca_5/1729 - \ca_6/133 $

$r=3$: \rule{0mm}{5mm} $ 17\sa_1/81252 +135\sa_2/117364 + 1775\sa_3/1056276 + 32\sa_4/20313$\\ \indent \rule{0mm}{5mm}$ + 167\sa_5/88023 + 203\sa_6/88023$

$r=4$: \rule{0mm}{5mm} The result is rather complicated and, therefore, omitted. For instance, the coefficient of $\ca_1$ is $- 186961973/560088713912$.

$n=14$: \rule{0mm}{5mm}

$r=1$: \rule{0mm}{5mm}$-\sa_4/7 - \sa_6/7$

$r=2$: \rule{0mm}{5mm} $13\ca_2/252  - \ca_4/252 + \ca_6/36$

$r=3$: \rule{0mm}{5mm} $3\sa_2/2336 + 43\sa_4/16352 + 51\sa_6/16352$

$r=4$: \rule{0mm}{5mm} $-11839\ca_2/25880400  + 2239\ca_4/25880400- 7489\ca_6/25880400$

$n=15$: \rule{0mm}{5mm}

$r=1$: \rule{0mm}{5mm} $- \sa_1/6 -\sa_2/10 - \sa_4/30 - \sa_7/10$

$r=2$: \rule{0mm}{5mm}$ \ca_1/192+7\ca_2/192 + 11\ca_4/960 - 19\ca_7/960 $

$r=3$: \rule{0mm}{5mm}$ 2797\sa_1/981120 +377\sa_2/327040 + 899\sa_4/981120  + 937\sa_7/327040$

$r=4$: \rule{0mm}{5mm}$ 499\ca_1/79412736-108943\ca_2/397063680 - 56287\ca_4/397063680$\\ \indent \rule{0mm}{5mm} $+ 14179\ca_7/79412736 $

%%%%%%%%%%%%%%%%%%%%%%%%%%%%%%%%%%%%%%%%%%%%%
\section*{5. Series giving rational values}
%%%%%%%%%%%%%%%%%%%%%%%%%%%%%%%%%%%%%%%%%%%%%

If $n$ is square-free, formula (\ref{3.26}) implies
\BE
\label{5.2}
  \CTRH\cdot S=(\BRTH)^T \:\MB{ and }\: \CTRH\cdot C=(\BRTH)^T
\EE
in the respective cases. We also use the connection (\ref{1.4}) between the numbers $\ctrh_j$ and $\DR_j$
and the representation (\ref{1.2}) of $\DR_j$ as an infinite sum. Then (\ref{5.2}) gives infinite series for the entries of the matrix $\BRTH$. We note the following cases.

\begin{prop} % Proposition 4 %%%%%%%%%%%%%%%%%%%%%%%%%%%%%%%%%%%%%%%%%%
\label{p4}

If $n\ge 3$ is square-free, then
\BE
\label{5.4}
(-1)^{(r-1)/2}\frac{2\pi^r}{n^r(r-1)!}\sum_{\btop{m>0,}{(m,n)=1}}\mu(m)\frac{\sin(2\pi m^*/n)}{m^r}= \BRTH_{1,1}
\EE
for odd numbers $r$, and
\BE
\label{5.6}
  (-1)^{r/2-1}\frac{2\pi^r}{n^r(r-1)!}\sum_{\btop{m>0,}{(m,n)=1}}\mu(m)\frac{\cos(2\pi m^*/n)}{m^r}=\BRTH_{1,1}.
\EE
for even numbers $r$.
\end{prop} %%%%%%%%%%%%%%%%%%%%%%%%%%%%%%%%%%%%%%%%%%%%%%%%%%%%%%%%%%%%%%%%%%

\MN
In the case $r=1$, we investigated this sum in \cite{Gi1}, however, in a way unrelated to the matrix $\BRTH$.

\MN
{\em Example.} For $n=35$, $r=3$, we have
\BD
   \BRTH_{1,1}=-\frac{4347647145233163511}{36746725032952512514560}=-0.0001183138671904\LD.
\ED
The right hand side of (\ref{5.4}), evaluated for $m\le 10000$, gives $-0.0001183138672025\LD$.

%%%%%%%%%%%%%%%%%%%%%%%%%%%%%%%%%%%%%%%%%%%%%%%%%%%%%%%%%%%%%%%%%%%%%
\section*{6. Eigenvalues and character coordinates}
%%%%%%%%%%%%%%%%%%%%%%%%%%%%%%%%%%%%%%%%%%%%%%%%%%%%%%%%%%%%%%%%%%%%%

In this section let $n$ be an integer, $n\ge 3$.
Let $a\in\Q(\zeta_n)$ and $\chi$ a Dirichlet character mod $n$ with conductor $f_{\chi}$. Let $\chi_f$ denote the Dirichlet character mod $f_{\chi}$ attached to $\chi$ (so $\chi_f(j)=\chi(j)$
if $(j,n)=1$). The $\chi$-{\em coordinate} $y(\chi|a)$ was introduced by Leopoldt (see \cite{Leo}) and is defined by the equation
\BD
  y(\chi|a)\tau(\OV{\chi}_f)=\sum_{\btop{1\le j\le n}{(j,n)=1}}\OV{\chi}(j)\sigma_j(a);
\ED
here $\OV{\chi}=\chi^{-1}$ is the complex-conjugate (or inverse) character of $\chi$ and $\tau(\chi_f)$ the {\em Gauss sum}
\BD
  \sum_{j=1}^{f_\chi}\chi_f(j)\zeta_{f_{\chi}}^{-j},
\ED
see \cite{Gi2}. Character coordinates are important for the study of Galois modules. For instance, the $\Q$-dimension of
\BD
   \sum_{\btop{1\le j\le n}{(j,n)=1}}\Q\sigma_j(a)
\ED
equals the number of characters (mod $n$) such that $y(\chi|a)\ne 0$.

An important property of the $\chi$-coordinate says
\BE
\label{6.2}
y(\chi|\sigma_k(a))=\chi(k)y(\chi|a).
\EE
We also need the {\em reconstruction formula}
\BE
\label{6.4}
  a=\frac 1{\phi(n)}\sum_{\chi \in\XX}y(\chi|a)\tau(\OV{\chi}_f),
\EE
where $\XX$ is the group of Dirichlet characters mod $n$.
This formula easily follows from the orthogonality relation for Dirichlet characters.

For a positive integer $r$, put
\BD
 \XXR=\{\chi\in\XX;\chi(-1)=(-1)^r\}.
\ED
Therefore, $\XXR$ is the set of odd (even, respectively) characters if $r$ is odd (even, respectively).
Suppose, henceforth, that $a\in\Q(\zeta_n)$ is purely imaginary if $r$ is odd, and real if $r$ is even. 
Then $y(\chi|a)\ne 0$ only if $\chi\in\XXR$.

We consider the matrix
\BE
\label{6.6}
  A=(\sigma_{jk^*}(a))_{j,k\in\RR}.
\EE
Further, we suppose that the set $\XXR$ is arranged in some order, so that we can  form the matrix
\BD
  X=\sqrt{2/\phi(n)}(\chi(k))_{k\in\RR,\chi\in\XXR}.
\ED
This matrix is {\em unitary}, i.e., $X\OV{X}^T=I$, where $I$ denotes the unit matrix (see \cite{Gi4}).
Let $\Delta$ be the diagonal matrix defined by
\BD
  \Delta=\left(\frac 12y(\chi|a)\tau(\OV{\chi}_f)\delta_{\chi,\psi}\right)_{\chi,\psi\in\XXR}.
\ED
Then the reconstruction formula (\ref{6.4}) and formula (\ref{6.2}) show
\BE
\label{6.8}
  X\Delta\OV{X}^T=A.
\EE
In other words, the character coordinates of $a$ are essentially the eigenvalues of $A$.

In the special case $a=\ctr_1$ we have $A=\CTR$, a matrix which, by (\ref{6.8}), is unitarly congruent to
\BD
 \Gamma=\left(\frac 12y(\chi|\ctr_1)\tau(\OV{\chi}_f)\delta_{\chi,\psi}\right)_{\chi,\psi\in\XXR}.
\ED
Since $\CTR$ is invertible, its inverse $\CTRH$ is unitarily congruent to $\Gamma^{-1}$.

Note, however, that $\CTRH$ is {\em not} of the form $(\sigma_{jk^*}(\ctrh_1))_{j,k\in\RR}$ (i.e., it does not go with (\ref{6.6})).
But $(\CTRH)^T$ has this form. From the identities
\BD
  \CTRH=X\Gamma^{-1}\OV X^T, \enspace (\CTRH)^T=X\WH{\Gamma}\OV X^T,
\ED
where $\WH{\Gamma}$ is the diagonal matrix
\BE
\label{6.10}
 \WH{\Gamma}=\left(\frac 12y(\chi|\ctrh_1)\tau(\OV{\chi}_f)\delta_{\chi,\psi}\right)_{\chi,\psi\in\XXR}
\EE
we obtain
\BD
 \WH{\Gamma}=\OV X^T \OV X\Gamma^{-1}X^TX.
\ED
Some calculation shows
\BD
  \WH{\Gamma}=(2y(\OV{\chi}|\ctr_1)^{-1}\tau(\chi_f)^{-1}\delta_{\chi,\psi})_{\chi,\psi\in\RR}
\ED
and, in view of (\ref{6.10}),
\BD
  y(\chi|\ctrh_1)\tau(\OV{\chi}_f)=4y(\OV{\chi}|\ctr_1)^{-1}\tau(\chi_f)^{-1}.
\ED
Moreover, we use the identity
\BD
 \tau(\chi_f)\tau(\OV{\chi}_f)=(-1)^rf_{\chi},
\ED
see \cite[p. 36]{Wa}. Then we may write
\BE
\label{6.12}
   y(\chi|\ctrh_1)=\frac{4(-1)^r}{f_{\chi}}y(\OV{\chi}|\ctr_1)^{-1}.
\EE
In other words, the character coordinates of the number $\ctrh_1$ are, in the main, the inverses of the character coordinates of the cotangent number $\ctr_1$.
We think that also for this reason we may speak of \BQ inverse\EQ cotangent numbers.

The character coordinates of $\ctr_1$ are given in \cite{Gi2}. Indeed,
\BE
\label{6.14}
  y(\chi|\ctr_1)=\left\{
                   \begin{array}{ll}
                     0, & \hbox{if $\chi\not\in\XXR$;} \\
                     (2n/f_{\chi})^r\prod_{p\DIV n}(1-\OV{\chi}_f(p)/p^r)\BRCHIF/r, & \hbox{if $\chi\in\XXR$.}
                   \end{array}
                 \right.
\EE
Here $\BRCHIF$ is the generalized Bernoulli number
\BD
  \BRCHIF=f_{\chi}^{r-1}\sum_{j=1}^{f_{\chi}}\chi_f(j)\BR(j/f_{\chi}),
\ED
where $\BR$ is, as above, the $r$th Bernoulli polynomial.
Now the reconstruction formula (\ref{6.4}), combined with (\ref{6.12}) and (\ref{6.14}), gives
the following theorem.

\begin{theorem} % Theorem 2 %%%%%%%%%%%%%%%%%%%%%%%%%%%%%%%%%%%%%%%%%%%%%%%%%%
\label{t2}

For all integers $n\ge 3$ and $r\ge 1$,
\BD
 \ctrh_1=(-1)^r\frac{2^{2-r}r}{\phi(n)n^r}\sum_{\chi\in\XXR}f_{\chi}^{r-1}\prod_{p\DIV n}(1-\chi_f(p)/p^r)^{-1}\frac{\tau(\OV{\chi}_f)}{\BRCHIFB}.
\ED

\end{theorem} %%%%%%%%%%%%%%%%%%%%%%%%%%%%%%%%%%%%%%%%%%%%%%%%%%%%%%%%%%%%%%%%

\MN
{\em Example.}
In the case $n=15$, $r=3$ there are four odd Dirichlet characters $\chi_1$,\LD,$\chi_4$, mod $n$ of conductors $3$, $5$, $5$, $15$, respectively. Then
\BEA
  \ct3h_1&=&-\frac{3}{16\cdot 15^3}\left(9(1-\chi_{1,3}(5)/125)^{-1}\frac{\tau(\chi_{1,3})}{\B3_{\chi_{1,3}}}+25(1-\chi_{2,5}(3)/27)^{-1}\frac{\tau(\chi_{3,5})}{\B3_{\chi_{3,5}}}\right.\\
         && \left.+25(1-\chi_{3,5}(3)/27)^{-1}\frac{\tau(\chi_{2,5})}{\B3_{\chi_{2,5}}} +225\frac{\tau(\chi_{4,15})}{\B3_{\chi_{4,15}}}\right).
\EEA
The second subscript of the $\chi$'s is always the conductor.
So $\chi_{1,3}$ is the character mod 3 attached to $\chi_1$. It is given by $\chi_{1,3}(2)=-1$. In the same way $\chi_{2,5}$ is given by $\chi_{2,5}(2)=i$,
and $\chi_{3,5}=\OV{\chi_{2,5}}$. Finally, $\chi_{4,15}$ is given by $\chi_{4,15}(2)=1$ and $\chi_{4,15} (-1)=-1$. The values of the Gauss sums are as follows:
$\tau(\chi_{1,3})=-s_{1,3}$, $\tau(\chi_{2,5})=-s_{1,5}-is_{2,5}$, $\tau(\chi_{3,5})=-s_{1,5}+is_{2,5}$, and $\tau(\chi_{4,15})=-s_{1,15}-s_{2,15}-s_{4,15}+s_{7,15}$.
Here $s_{j,m}$ stands for $2i\sin(2\pi j/m)$. The corresponding Bernoulli numbers take the following values: $\B3_{\chi_{1,3}}=2/3$, $\B3_{\chi_{2,5}}=(12+6i)/5$,
$\B3_{\chi_{3,5}}=(12-6i)/5$, and $\B3_{\chi_{4,15}}=48$. Altogether, we obtain

$\ct3h_1=i(\sin(\pi/3)/672 + 11\sin(2\pi/5)/5840 + \sin(\pi/5)/1168 + \sin(2\pi/15)/1920$ \rule{0mm}{5mm}\\
\indent        $+ \sin(4\pi/15)/1920+ \sin(7\pi/15)/1920 -  \sin(\pi/15)/1920)$. \rule{0mm}{5mm}

\MN
So we have expressed $\ct3h_1$ as a rational linear combination of powers of $\zeta_{15}$, but these powers are not $\Q$-linearly independent here -- in contrast with the respective result of Section 4.

\MN
{\em Remarks}.  1. Let $n=p\ge 3$ be a prime and $r$ odd. Then all odd Dirichlet characters mod $p$ have the conductor $p$. Accordingly, Theorem \ref{t2} says
\BD
 \ctrh_1=-\frac{2^{2-r}}{(p-1)p}\sum_{\chi\in\XXR}\frac{\tau({\chi}_f)}{\BR_{\chi_f}}.
\ED
This formula looks remarkably simple. Maybe Hecke envisioned something of this kind when he said that the case of a prime number $n$ is easier than the composite case (see Section 1).

2. The character coordinates of the numbers $\ctr_1$ are closely related to those of $i^r\cot(\pi/n)^r$, since $\cot(x)^r$ is a rational linear combination of derivatives of $\cot(x)$.
The character coordinates of the last-mentioned numbers can be found in \cite[Cor. 4.4]{Is} and, for primitive characters, in \cite[Cor. 2.19]{Fr}, but without any reference to this name and, consequently,
to our paper \cite{Gi2}.

\bigskip
\centerline{\bf Conflicting interests}

\MN
The author declares that there are no conflicting interests.

%%%%%%%%%%%%%%%%%%%%%%%%%%%%%%%%%%%%%%%%%%%%%%%%%%%%%
%%%%%%%%%%%%%%%%%%%%%%%%%%%%%%%%%%%%%%%%%%%%%%%%%%%%%%%%%%%%%%%%%%%%%%%%%%

\MN
Kurt Girstmair\\
Institut f\"ur Mathematik \\
Universit\"at Innsbruck   \\
Technikerstr. 13/7        \\
A-6020 Innsbruck, Austria \\
Kurt.Girstmair@uibk.ac.at

\end{document}